\documentclass[11pt,reqno]{amsart}
\usepackage[T1]{fontenc}
\usepackage{lmodern}
\usepackage{amsmath,amssymb,amsthm,mathtools}
\usepackage{microtype}
\usepackage[a4paper,margin=30mm]{geometry}
\usepackage{enumitem}
\usepackage{url}
\usepackage[hidelinks,pdfencoding=auto,psdextra]{hyperref}
\hypersetup{pdftitle={Failure of bounded dbar-solvability on pseudoconvex domains of finite type},pdfauthor={Yun-Heng Du},pdfkeywords={dbar equation, sup-norm estimates, finite-type domains, singular complex curves}}
\numberwithin{equation}{section}
\newtheorem{theorem}{Theorem}[section]
\newtheorem{proposition}[theorem]{Proposition}

\newcommand{\C}{\mathbb C}
\newcommand{\R}{\mathbb R}

\newcommand{\db}{\bar\partial}
\DeclareMathOperator{\dist}{dist}
\setlist[enumerate,1]{label=\textup{(\roman*)},leftmargin=2.3em}
\title[Failure of bounded $\bar\partial$-solvability on pseudoconvex domains of finite type]{Failure of bounded $\bar\partial$-solvability\\on pseudoconvex domains of finite type}
\author{Yun-Heng Du}
\address[]{Academy of Mathematics and Systems Science, Chinese Academy of Sciences, Beijing 100190, China}
\email{duyunheng22@mails.ucas.ac.cn}
\date{}
\begin{document}
\raggedbottom
\begin{abstract}
We show that for every integer $n\ge3$ there is a bounded pseudoconvex domain $\Omega\subset\C^n$ with smooth real-analytic boundary whose D'Angelo $1$-type has maximum $16$ and a $\bar\partial$-closed $(0,1)$-form $f$, continuous on $\overline\Omega$ and smooth in $\Omega$, such that $\bar\partial g=f$ has no solution in $L^\infty(\Omega)$. On each of these domains, we also construct a uniformly bounded sequence of smooth $\db$-closed $(0,1)$-forms for which the infimum of the $L^\infty$ norms of bounded solutions tends to infinity. Each form in this sequence, together with one of its solutions, extends smoothly to a neighborhood of the closure. The construction uses a singular curve both to constrain holomorphic corrections and to separate the scales of the cutoff error and the resulting lower bound.
\end{abstract}
\maketitle
\section{Introduction}
Bounded solutions of the $\db$-equation allow errors in local holomorphic constructions to be corrected with uniform control. For a bounded pseudoconvex domain $D\Subset\C^n$, the sup-norm estimate requires a constant $C_D$ such that every bounded, distributionally $\db$-closed $(0,1)$-form $f$ admits a solution $g$ of
\begin{equation}\label{eq:estimate}
\db g=f,\qquad \|g\|_{L^\infty(D)}\le C_D\|f\|_{L^\infty(D)}.
\end{equation}
The classical theory establishes this estimate on strongly pseudoconvex domains; see \cite{Henkin,Ramirez,GrauertLieb,Kerzman}. The question addressed here is whether finite D'Angelo type suffices in dimensions $n\ge3$.

Finite D'Angelo type is closely tied to subellipticity of the $\db$-Neumann problem~\cite{Kohn,DAngelo,Catlin}. Subelliptic estimates control Sobolev regularity, whereas \eqref{eq:estimate} asks for control of a solution solely by the $L^\infty$ norm of the right-hand side. For bounded pseudoconvex domains with smooth real-analytic boundary, Diederich--Forn\ae ss~\cite{DFRealAnalytic} had already established subelliptic estimates in 1978. In complex dimension two, Forn\ae ss~\cite{Fornaess} obtained sup-norm estimates on a class containing the Kohn--Nirenberg domain. The sup-norm and H\"older estimates of Fefferman--Kohn~\cite{FeffermanKohn} establish \eqref{eq:estimate} on smooth bounded pseudoconvex domains of finite type; Range~\cite{Range} subsequently obtained H\"older estimates by integral-kernel methods. For bounded convex domains in $\C^2$ with $C^1$ boundary and without a finite-type assumption, Forn\ae ss--Kim~\cite[Problem~3]{FornaessKim} asked whether sup-norm estimates hold.

In higher dimensions, sup-norm estimates hold on ellipsoids~\cite{DFW} and convex finite-type domains~\cite{DFF,Cumenge,Yao}. The latter theory draws on the geometry studied by McNeal~\cite{McNeal} and the holomorphic support functions constructed by Diederich--Forn\ae ss~\cite{DFSupport}. Further estimates hold on lineally convex finite-type domains~\cite{DFischer}; local estimates are also available under a smooth diagonalizability assumption on the Levi form~\cite{FKM}. Gavosto~\cite{Gavosto} obtained local estimates on certain finite-type domains related to strongly pseudoconvex domains by proper mappings branched in a single coordinate.

Pseudoconvexity alone does not imply bounded solvability. Sibony~\cite{Sibony} constructed a smoothly bounded pseudoconvex domain in $\C^3$ and a $\db$-closed form continuous on its closure with no bounded solution. The boundary of his domain contains analytic discs, and he asked about the real-analytic boundary case~\cite[p.~236]{Sibony}. Berndtsson~\cite{Berndtsson} later obtained a counterexample on a smoothly bounded Hartogs domain in $\C^2$ by relating bounded solvability to weighted one-variable estimates. Related failures of $L^p$ estimates were obtained by Forn\ae ss--Sibony~\cite{FS}. 

Parrini--Selvaggi Primicerio~\cite[Theorem~3.1]{PSP} asserted a sup-norm estimate for finite covers satisfying a quantitative pseudoconvex-neighborhood condition. Barletta--Parrini~\cite[Theorem~3.1]{BP} asserted the existence of bounded solutions on smooth pseudoconvex finite covers of strongly pseudoconvex domains, with the covering map smooth up to the boundary, for forms in a weighted class that includes bounded forms smooth in the domain. Grundmeier--Simon--Stens\o nes~\cite[Main Theorem~1]{GSS} asserted a sup-norm estimate, together with H\"older estimates of positive order, on every bounded pseudoconvex domain in $\C^3$ with real-analytic boundary of finite D'Angelo type. These three bounded-solvability assertions are false, as the counterexamples below show.

We denote the vanishing order of a holomorphic germ $a(t)$ at $0$ by $\nu(a)$, with $\nu(0)=+\infty$. The order of a vector-valued holomorphic germ is the least order of its components. For a real-analytic germ $b(t,\bar t)=\sum_{j,k\ge0}b_{jk}t^j\bar t^k$ at $0$, we use
\[
\nu(b)=\min\{j+k:b_{jk}\ne0\}.
\] The D'Angelo $1$-type of a smooth real-analytic hypersurface $M$ at $p\in M$ is
\begin{equation}\label{eq:type}
\Delta_1(M,p)=\sup_\gamma\frac{\nu(\rho\circ\gamma)}{\nu(\gamma-p)},
\end{equation}
where $\rho$ is a real-analytic local defining function for $M$ and the supremum is over nonconstant holomorphic curve germs $\gamma:(\C,0)\to(\C^n,p)$.

We consider the domains
\begin{equation}\label{eq:domain}
\Omega_n\coloneqq\left\{(u,v,w,\zeta)\in\C^n:|w|^2+|u^4-v^5|^2+|v|^{16}+\sum_{\ell=1}^{n-3}|\zeta_\ell|^2<1\right\},\quad n\ge3.
\end{equation}
Here $\zeta=(\zeta_1,\ldots,\zeta_{n-3})$, and the sum is empty when $n=3$. Throughout the paper, forms are measured in the Euclidean coefficient norm, and $L^\infty$ norms are essential suprema with respect to Euclidean volume. Functions in $C^\infty(\overline\Omega_n)$ are restrictions of functions smooth on an open neighborhood of $\overline\Omega_n$. The spaces $C^0_{(0,1)}(X)$ and $C^\infty_{(0,1)}(X)$ consist of $(0,1)$-forms with coefficients in $C^0(X)$ and $C^\infty(X)$ respectively. We write
\[
L^\infty(\Omega_n)=\bigl\{g:\Omega_n\to\C\text{ measurable}:\|g\|_{L^\infty(\Omega_n)}=\operatorname*{ess\,sup}_{x\in\Omega_n}|g(x)|<\infty\bigr\}.
\]
Here $\operatorname*{ess\,sup}_{x\in\Omega_n}|g(x)|$ is the least $M\ge0$ such that $|g(x)|\le M$ for almost every $x\in\Omega_n$, and functions agreeing almost everywhere are identified.
\begin{theorem}\label{thm:main}
For every integer $n\ge3$, the domain $\Omega_n$ in \eqref{eq:domain} has the following properties.
\begin{enumerate}
\item It is bounded, connected, and pseudoconvex, with smooth real-analytic boundary. Its D'Angelo $1$-type is at most $16$ at every boundary point and equals $16$ at $(0,0,1,0)$.
\item With $k_0\coloneqq2^{134}$, there are constants $C_0,c>0$ and $\db$-closed forms $f_k\in C^\infty_{(0,1)}(\overline\Omega_n)$ for integers $k\ge k_0$ such that $\|f_k\|_{L^\infty(\Omega_n)}\le C_0$, and each bounded distributional solution $g\in L^\infty(\Omega_n)$ of $\db g=f_k$ satisfies
\begin{equation}\label{eq:mainlower}
\|g\|_{L^\infty(\Omega_n)}\ge c k^{1/320}.
\end{equation}
The forms $f_k$ also admit solutions in $C^\infty(\overline\Omega_n)$.
\item There is a distributionally $\db$-closed form $f\in C^0_{(0,1)}(\overline\Omega_n)\cap C^\infty_{(0,1)}(\Omega_n)$ for which $\db g=f$ has no solution in $L^\infty(\Omega_n)$.
\end{enumerate}
\end{theorem}
Theorem~\ref{thm:main} gives a negative answer to Sibony's question about real-analytic boundaries~\cite[p.~236]{Sibony} in every dimension $n\ge3$, with a fixed bound on the D'Angelo $1$-type. The domain $\Omega_n$ is a finite proper preimage of the unit ball and satisfies the quantitative pseudoconvex-neighborhood condition discussed above; these properties are verified in Section~\ref{sec:construction} and Appendix~\ref{app:covers}. Theorem~\ref{thm:main}(ii) disproves the uniform estimates asserted in \cite{PSP,GSS}, while Theorem~\ref{thm:main}(iii) disproves the bounded-solvability assertion in \cite{BP}.

Following Sibony's strategy of cutting off a singular function and obstructing holomorphic corrections~\cite[pp.~239--240]{Sibony}, we use a singular curve both to restrict those corrections and to separate geometric scales. On the curve $(u,v)=(t^5,t^4)$, restrictions of holomorphic functions have only powers $t^{5a+4b}$ with nonnegative integers $a,b$. The missing degree $11$ suggests the quotient $u^3/v$, whose restriction is $t^{11}$ and therefore cannot be canceled by a holomorphic correction. Since $u^4-v^5$ vanishes on this curve, slices of $\Omega_n$ with $w$ fixed near $1$ and $\zeta=0$ extend farther along the curve than in transverse directions, while the term $|v|^{16}$ keeps the boundary of finite type. This construction uses two complex tangential directions at $(0,0,1,0)$, whereas a hypersurface in $\C^2$ has only one.

We cut off $u^3/v$ near $v=0$ and multiply by $w^k$. Near $|w|=1$, the defining inequality makes $u^3$ small on the cutoff annulus; away from that region, $w^k$ suppresses the cutoff error. The resulting $\db$ derivatives remain uniformly bounded, while the surviving coefficient measured on larger circles along the singular curve grows. Decay of $w^k$ on compact subsets also allows the weighted series to remain smooth in the interior. The passage to a single form with no bounded solution uses Berndtsson's idea of summing a suitably chosen subsequence~\cite[Proposition~2.1]{Berndtsson}. Here the summands contribute nonnegatively to the selected Fourier coefficient, so the obstruction survives summation.

Section~\ref{sec:construction} establishes the geometry of $\Omega_n$ and constructs the sequence of bounded forms that rules out uniform sup-norm estimates. A weighted sum of a subsequence yields the counterexample to bounded solvability in Section~\ref{sec:single}. In Appendix~\ref{app:covers}, we verify the hypotheses of the earlier bounded-solvability assertions and identify the failures in their kernel and division estimates.

\section{Unbounded solution norms for smooth forms}\label{sec:construction}

Fix $n\ge3$ and consider
\[
\Omega_n=\{(u,v,w,\zeta)\in\C^n:|w|^2+|u^4-v^5|^2+|v|^{16}+|\zeta|^2<1\},
\]
where $\zeta=(\zeta_1,\ldots,\zeta_{n-3})$ and $|\zeta|^2=\sum_{\ell=1}^{n-3}|\zeta_\ell|^2$. For $n=3$, the $\zeta$ coordinates are absent and $|\zeta|^2=0$. Write $x\coloneqq(u,v,w,\zeta)$. The polynomial map
\begin{equation}\label{eq:map} F(x)\coloneqq(u^4-v^5,v^8,w,\zeta)\end{equation}
gives the representation
\begin{equation}\label{eq:rho} Q(x)\coloneqq|F(x)|^2,\qquad \rho(x)\coloneqq Q(x)-1,\qquad \Omega_n=\{\rho<0\}.\end{equation}
\begin{proposition}\label{prop:geometry}
The domain $\Omega_n$ has the following properties.
\begin{enumerate}
\item It is bounded, connected, and pseudoconvex, and its boundary is smooth and real analytic.
\item At every $p\in\partial\Omega_n$, the D'Angelo $1$-type satisfies $\Delta_1(\partial\Omega_n,p)\le16$, with equality at $p=(0,0,1,0)$.
\item The map $F:\C^n\to\C^n$ in \eqref{eq:map} is proper and surjective. Every $y=(y_1,y_2,y_3,y')$ satisfying $y_2\ne0$ and $y_1^8\ne y_2^5$ has exactly $32$ distinct preimages under $F$. The restriction of $F$ maps $\Omega_n$ properly onto the unit ball $B^n$.
\end{enumerate}
\end{proposition}
\begin{proof}
(i) The inequality $Q<1$ gives $|v|<1$, $|w|<1$, and $|\zeta|<1$, while $|u|^4\le |u^4-v^5|+|v|^5<2$. Thus $\Omega_n$ is bounded. For $0\le s\le1$, the maps $\delta_s(u,v,w,\zeta)\coloneqq(s^5u,s^4v,sw,s\zeta)$ satisfy
\begin{equation}\label{eq:contraction}
Q(\delta_sx)=s^{40}|u^4-v^5|^2+s^{64}|v|^{16}+s^2(|w|^2+|\zeta|^2)\le s^2Q(x).
\end{equation}
They contract $\Omega_n$ to $0$ through maps taking values in $\Omega_n$, so $\Omega_n$ is connected. At a point with $Q(x)=1$,
\[
\left.\frac{d}{ds}Q(\delta_sx)\right|_{s=1}=40|u^4-v^5|^2+64|v|^{16}+2|w|^2+2|\zeta|^2\ge2.
\]
Consequently $d\rho\ne0$ on $\{Q=1\}$, so $\partial\Omega_n=\{Q=1\}$ is a smooth real-analytic hypersurface. Since $Q$ is a sum of squared absolute values of holomorphic polynomials, $\rho$ is plurisubharmonic. The function $-\log(-\rho)$ is a plurisubharmonic exhaustion of $\Omega_n$, which proves pseudoconvexity.

(ii) Consider a nonconstant holomorphic curve germ $\gamma(t)=(u(t),v(t),w(t),\zeta(t))$ with $\gamma(0)=p=(u_0,v_0,w_0,\zeta_0)\in\partial\Omega_n$. Set $A(t)\coloneqq F(\gamma(t))-F(p)$. We first prove
\begin{equation}\label{eq:orderA}\nu(A)\le8\nu(\gamma-p).\end{equation}
Write $\alpha\coloneqq\nu(u-u_0)$, $\beta\coloneqq\nu(v-v_0)$, $q\coloneqq u^4-v^5$, and $q_0\coloneqq u_0^4-v_0^5$. We claim that
\begin{equation}\label{eq:ordertwo}
\min\{\nu(q-q_0),\nu(v^8-v_0^8)\}\le8\min\{\alpha,\beta\}.
\end{equation}

Let $a\coloneqq\nu(u^4-u_0^4)$. Then $a=\alpha$ if $u_0\ne0$ and $a=4\alpha$ otherwise. If $v_0\ne0$, both $v^5-v_0^5$ and $v^8-v_0^8$ have order $\beta$, so
\[
\min\{\nu(q-q_0),\nu(v^8-v_0^8)\}
=\min\{a,\beta\}\le4\min\{\alpha,\beta\}.
\]
If $v_0=0$ and $\beta\le\alpha$, the component $v^8$ has order $8\beta$. If $v_0=0$ and $\alpha<\beta$, then $a\le4\alpha<5\beta=\nu(v^5)$, so $\nu(q-q_0)=a\le4\alpha$. This proves \eqref{eq:ordertwo}.

Denote by $b$ the least vanishing order among the remaining components $w-w_0$ and $\zeta_\ell-\zeta_{\ell,0}$. These are also components of $\gamma-p$, so $\nu(\gamma-p)=\min\{\alpha,\beta,b\}$. By \eqref{eq:ordertwo},
\[
\nu(A)\le\min\{8\min\{\alpha,\beta\},b\}
\le8\min\{\alpha,\beta,b\}=8\nu(\gamma-p),
\]
which proves \eqref{eq:orderA}. In particular, $m\coloneqq\nu(A)$ is finite. Because $|F(p)|^2=1$,
\begin{equation}\label{eq:mixed}
\rho(\gamma(t))=2\operatorname{Re}\langle A(t),F(p)\rangle+|A(t)|^2.
\end{equation}
Write $A(t)=ct^m+O(t^{m+1})$ with $c\in\C^n\setminus\{0\}$. The coefficient of $t^m\bar t^m$ in \eqref{eq:mixed} is $|c|^2>0$, since $2\operatorname{Re}\langle A(t),F(p)\rangle$ contains no mixed monomials. Therefore $\nu(\rho\circ\gamma)\le2m\le16\nu(\gamma-p)$. At $p_*\coloneqq(0,0,1,0)$, the curve $\gamma_*(t)\coloneqq(t^5,t^4,1,0)$ satisfies $\nu(\gamma_*-p_*)=4$ and $\rho\circ\gamma_*=|t|^{64}$. Thus $\Delta_1(\partial\Omega_n,p_*)=16$.

(iii) For $y=(y_1,y_2,y_3,y')\in\C^n$, the equations $F(x)=y$ are
\[v^8=y_2,\qquad u^4=y_1+v^5,\qquad w=y_3,\qquad\zeta=y'.\]
This system is solvable for every $y$. If $y_2\ne0$ and $y_1^8\ne y_2^5$, there are eight distinct roots for $v$ and four distinct roots for $u$ for each $v$. These inequalities describe an open dense set of points $y$, each with exactly $32$ preimages. If $y$ ranges over a bounded set, the equations bound $v$, then $u$, and also $w$ and $\zeta$. The inverse image of a compact set is thus closed and bounded, so $F$ is proper and restricts to a proper surjective map $\Omega_n=F^{-1}(B^n)\to B^n$.
\end{proof}

Proposition~\ref{prop:geometry} establishes Theorem~\ref{thm:main}(i). To construct the forms in Theorem~\ref{thm:main}(ii), choose a real-valued function $\chi\in C^\infty(\R)$ such that
\begin{equation}\label{eq:cutoff}
0\le\chi\le1,\qquad\chi(s)=0\ (s\le1),\qquad\chi(s)=1\ (s\ge4).
\end{equation}
For every integer $k\ge2$, consider the smooth function
\begin{equation}\label{eq:Vk}
V_k(u,v,w,\zeta)\coloneqq k^{7/40}w^k\frac{u^3}{v}\chi(k^{1/5}|v|^2).
\end{equation}
Differentiation gives the following smooth $\db$-closed forms on $\C^n$:
\begin{equation}\label{eq:fk}
f_k\coloneqq\db V_k=k^{3/8}w^ku^3\chi'(k^{1/5}|v|^2)\,d\bar v.
\end{equation}

\begin{proposition}\label{prop:data}
The forms in \eqref{eq:fk} satisfy
\begin{equation}\label{eq:databound}
\|f_k\|_{L^\infty(\Omega_n)}\le C_\chi,\qquad C_\chi\coloneqq(32^{3/4}+1)\sup_{s\in\R}|\chi'(s)|,\qquad k\ge2.
\end{equation}
\end{proposition}
\begin{proof}
Since $\chi'$ vanishes outside $[1,4]$, the form $f_k$ vanishes unless $v$ belongs to the annulus $\{v\in\C:k^{-1/10}\le|v|\le2k^{-1/10}\}$. For $x=(u,v,w,\zeta)\in\Omega_n$ with $v$ in this annulus, one has
\[|u^4-v^5|\le(1-|w|^2)^{1/2},\qquad |u|^4\le32k^{-1/2}+(1-|w|^2)^{1/2}.\]
Taking the $3/4$ power gives
\begin{equation}\label{eq:ubound}
|u|^3\le\left(32k^{-1/2}+(1-|w|^2)^{1/2}\right)^{3/4}\le32^{3/4}k^{-3/8}+(1-|w|^2)^{3/8}.
\end{equation}
Since $|w|^2<1$ and $y^{3/8}e^{-ky/2}\le(3/(4e))^{3/8}k^{-3/8}$ for $y\ge0$, we obtain that
\[
|w|^k(1-|w|^2)^{3/8}\le e^{-k(1-|w|^2)/2}(1-|w|^2)^{3/8}\le k^{-3/8}.
\]
Substituting this bound into \eqref{eq:ubound} and \eqref{eq:fk} gives
\[
|f_k|\le\sup_{s\in\R}|\chi'(s)|\left(32^{3/4}|w|^k+1\right)\le C_\chi,
\]
which proves \eqref{eq:databound}.
\end{proof}

For every integer $k\ge2$, consider the holomorphic map
\begin{equation}\label{eq:Psi}\Psi_k(t)\coloneqq(t^5,t^4,s_k,0),\qquad s_k\coloneqq\sqrt{1-1/k}.
\end{equation}
The image under $\Psi_k$ of the closed disc $\{|t|\le T_k\}$, where $T_k\coloneqq\tfrac12 k^{-1/64}$, is compactly contained in $\Omega_n$, since
\begin{equation}\label{eq:disc}
Q(\Psi_k(t))=s_k^2+|t|^{64}\le1-\frac1k+\frac{2^{-64}}k<1.
\end{equation}
On $\{|t|=T_k\}$, one has
\begin{equation}\label{eq:threshold}
k^{1/5}|t^4|^2=k^{1/5}T_k^8=2^{-8}k^{3/40}\ge4\qquad(k\ge2^{134}).
\end{equation}
Hence $\chi(k^{1/5}|t^4|^2)=1$ on this circle for $k\ge k_0$.
\begin{proof}[Proof of Theorem~\ref{thm:main}(ii)]
The forms $f_k=\db V_k$ satisfy $\|f_k\|_{L^\infty(\Omega_n)}\le C_\chi$ by Proposition~\ref{prop:data}, and each $V_k$ is smooth on $\C^n$. Fix an integer $k\ge k_0$ and a bounded distributional solution $g$ of $\db g=f_k$ on $\Omega_n$.
Since $\db(g-V_k)=0$ in distributions, also $\Delta(g-V_k)=0$. Weyl's lemma gives a smooth representative $H$ of $g-V_k$, which is holomorphic because $\db H=0$. Choosing the smooth representative $V_k+H$ of $g$, one has
\[
|g(x)|\le\|g\|_{L^\infty(\Omega_n)}\qquad(x\in\Omega_n).
\]

On $\{|t|=T_k\}$, \eqref{eq:Vk} and \eqref{eq:threshold} give
\[
V_k(\Psi_k(t))=k^{7/40}s_k^kt^{11}.
\]
By \eqref{eq:disc}, the function $t\mapsto H(\Psi_k(t))$ is holomorphic on a neighborhood of the closed disc $|t|\le T_k$. The Taylor expansion of $H(u,v,s_k,0)$ in $(u,v)$ contains only monomials $u^av^b$ with nonnegative integers $a,b$, which become $t^{5a+4b}$ after substitution. Since $5a+4b=11$ has no nonnegative integer solutions, the coefficient of $t^{11}$ is zero. Therefore,
\begin{equation}\label{eq:holomorphic-moment}
\frac1{2\pi}\int_0^{2\pi}H(\Psi_k(T_ke^{i\theta}))e^{-11i\theta}\,d\theta=0.
\end{equation}

Consequently,
\[
\begin{aligned}
\|g\|_{L^\infty(\Omega_n)}
&\ge\left|\frac1{2\pi}\int_0^{2\pi}g(\Psi_k(T_ke^{i\theta}))e^{-11i\theta}\,d\theta\right|\\
&=k^{7/40}s_k^kT_k^{11}
=2^{-11}(1-1/k)^{k/2}k^{1/320}
\ge e^{-1}2^{-11}k^{1/320}.
\end{aligned}
\]
This proves Theorem~\ref{thm:main}(ii), with $C_0=C_\chi$, $c=e^{-1}2^{-11}$, and $k_0=2^{134}$.
\end{proof}

\section{Nonexistence of bounded solutions}\label{sec:single}
Consider the form
\begin{equation}\label{eq:fsum}
f\coloneqq\sum_{j=1}^\infty2^{-j}f_{k_j},\qquad k_j\coloneqq2^{640j}\quad(j\ge1).
\end{equation}
By \eqref{eq:databound}, the series converges uniformly on $\overline\Omega_n$, so $f\in C^0_{(0,1)}(\overline\Omega_n)$ and $\|f\|_{L^\infty(\Omega_n)}\le C_\chi$.

\begin{proof}[Proof of Theorem~\ref{thm:main}(iii)]
Set
\begin{equation}\label{eq:Vsum}V\coloneqq\sum_{j=1}^\infty2^{-j}V_{k_j}.\end{equation}
For each compact set $K\Subset\Omega_n$, choose $q\in(0,1)$ with $|w|\le q$ on $K$. With $B(v)\coloneqq\chi(|v|^2)/v$, extended by zero near $v=0$, one has
\[
\frac{\chi(k^{1/5}|v|^2)}v=k^{1/10}B(k^{1/10}v).
\]
All real derivatives of $B$ are bounded. Hence every real derivative of order $a$ of $k^{1/10}B(k^{1/10}v)$ is bounded by $C_a k^{(a+1)/10}$, uniformly for $v\in\C$. It follows from \eqref{eq:Vk} that, for every multi-index $\alpha$,
\[
\sup_K|D^\alpha V_k|\le C_{K,\alpha}k^{A_\alpha}q^k
\]
for some $A_\alpha\ge0$. Applying the Weierstrass convergence criterion to $\sum_{j=1}^\infty2^{-j}D^\alpha V_{k_j}$ for each $\alpha$, one obtains a smooth sum $V$ in \eqref{eq:Vsum} with $f=\db V\in C^\infty_{(0,1)}(\Omega_n)$ and $\db f=0$.
Suppose that $g\in L^\infty(\Omega_n)$ satisfies $\db g=f$. As in the proof of Theorem~\ref{thm:main}(ii), choose the smooth representative $g=V+H$, where $H$ is holomorphic on $\Omega_n$, so that $|g(x)|\le M\coloneqq\|g\|_{L^\infty(\Omega_n)}$ for all $x\in\Omega_n$.

For a fixed $j$, \eqref{eq:disc} gives
\[
\Psi_{k_j}\bigl(\{|t|\le T_{k_j}\}\bigr)\Subset\Omega_n.
\]
Set
\begin{equation}\label{eq:Jj}
J_j\coloneqq\frac1{2\pi}\int_0^{2\pi}g\bigl(\Psi_{k_j}(T_{k_j}e^{i\theta})\bigr)e^{-11i\theta}\,d\theta.
\end{equation}
Then $|J_j|\le M$. Equation \eqref{eq:holomorphic-moment} shows that $H$ contributes nothing to $J_j$. The series \eqref{eq:Vsum} converges uniformly on $\Psi_{k_j}(\{|t|=T_{k_j}\})$, so termwise integration and \eqref{eq:threshold} give
\[
J_j=T_{k_j}^{11}\sum_{\ell=1}^\infty2^{-\ell}k_\ell^{7/40}s_{k_j}^{k_\ell}\chi(k_\ell^{1/5}T_{k_j}^8)
\ge2^{-j}k_j^{7/40}s_{k_j}^{k_j}T_{k_j}^{11}
\ge e^{-1}2^{-11}2^j.
\]
This contradicts $|J_j|\le M$ as $j\to\infty$.
\end{proof}

\appendix
\section{Finite-cover and holomorphic division estimates}\label{app:covers}
\subsection{Verification of the finite-cover hypotheses}
The domain $\Omega_n$ satisfies the geometric hypotheses of Parrini--Selvaggi Primicerio~\cite[Theorem~3.1]{PSP} and Barletta--Parrini~\cite[Theorem~3.1]{BP}. Indeed, Proposition~\ref{prop:geometry} gives a smooth finite proper preimage of the ball under the polynomial map $F$, of degree $32$. The neighborhood condition in \cite[Definition~1.1]{PSP} also holds: with $K\coloneqq\overline\Omega_n$ and $K_r\coloneqq\{x:\dist(x,K)<r\}$, there are pseudoconvex neighborhoods $D_j$ with intersection $K$ and a fixed $p\ge1$ such that, for every $0<\epsilon<1/2$, some $j$ satisfies
\begin{equation}\label{eq:neighborhood}
K_{\epsilon^p}\subset D_j\subset K_\epsilon.
\end{equation}
By Diederich--Forn\ae ss~\cite[Theorem~1]{DF}, $K$ is uniformly $H$-convex, since $F$ is a finite proper holomorphic map onto the ball that extends smoothly to a neighborhood of $\overline\Omega_n$. This property gives the quantitative neighborhood condition \eqref{eq:neighborhood}.

Consequently, the smooth forms in Theorem~\ref{thm:main}(ii) contradict the sup-norm estimate of \cite[Theorem~3.1]{PSP}; each $f_k$ has only one nonzero coefficient, so the maximum-coefficient norm used there coincides with our Euclidean norm. The form $f$ in Theorem~\ref{thm:main}(iii), which is smooth in the domain and continuous on its closure, satisfies the weighted hypothesis of \cite[Theorem~3.1]{BP}. Explicitly, their weighted condition is
\[
\sup_{\Omega_n}|\rho|^\alpha\bigl(|f|+|\rho|^{-1/(2M)}|f\wedge\db\rho|\bigr)<\infty,
\qquad \alpha<1/M.
\]
Here $M$ is the degree of the cover and $f$ denotes a $(0,1)$-form. For $M=32$, every bounded $f$ satisfies it with $\alpha=3/128$, since $\alpha-1/(2M)=1/128>0$ and $\db\rho$ is bounded. Thus Theorem~\ref{thm:main}(iii) contradicts their bounded-solvability assertion.

\subsection{The higher-dimensional gap in Parrini--Selvaggi Primicerio}
The gap in the proof of Parrini--Selvaggi Primicerio~\cite[Theorem~3.1]{PSP} concerns derivatives of the support vector in dimensions $n\ge3$. In their notation, $D$ is the domain, $g$ is the closed $(0,1)$-form, $\zeta$ is the integration variable, and $z$ is the evaluation point. In the construction of \cite[p.~467]{PSP}, write $p=(p_1,\ldots,p_n)$ and denote their scalar support function by $\mathcal P(\zeta,z)\coloneqq\sum_j(\zeta_j-z_j)p_j(\zeta,z)$. Set
\[
B\coloneqq\frac{\bar\zeta-\bar z}{|\zeta-z|^2},\qquad b\coloneqq\frac{p}{\mathcal P}.
\]
Their homotopy joins $B$ to $b$, and the kernel for solving the $(0,1)$-equation has one vector column and $n-1$ differential columns. In dimension two, integration in the homotopy parameter uses the only differential column, leaving no derivatives of $p$. The calculation in \cite[p.~472]{PSP} accordingly splits the solution into the Bochner--Martinelli volume term $J_1$ and a boundary transition term $J_2$, and gives
\[
|J_2(z)|\le c\|g\|\int_{\partial D}\frac{d\sigma(\zeta)}{|\mathcal P(\zeta,z)|\,|\zeta-z|}.
\]

In dimension three, one spatial differential column remains after the homotopy integration. The boundary transition form then contains terms proportional to $\det(B,b,\db_\zeta B)$ and $\det(B,b,\db_\zeta b)$. Since
\[
\db_\zeta b=\frac{\db_\zeta p}{\mathcal P}-\frac{p\,\db_\zeta\mathcal P}{\mathcal P^2},
\qquad
\det(B,b,\db_\zeta b)=\frac{\det(B,p,\db_\zeta p)}{\mathcal P^2},
\]
the latter term involves derivatives of $p$ and a second power of the support function in the denominator. The estimate for $J_2$ in \cite[p.~472]{PSP} does not control this term: an additional cancellation or boundary integral estimate is needed, but none is supplied. This is the missing step in extending their two-dimensional calculation to the general-dimensional assertion of Theorem~3.1. Our Theorem~\ref{thm:main}(ii) disproves that assertion for every $n\ge3$.

\subsection{Failure of the Barletta--Parrini kernel decomposition}
In Barletta--Parrini~\cite{BP}, the specific error occurs in the passage from the kernel in \cite[Proposition~3.1]{BP} to the bounded-coefficient decomposition asserted in \cite[p.~125, (3.3)]{BP}. We examine this passage on $\Omega_3$. Write $x\coloneqq(u,v,w)$, $y\coloneqq(u',v',w')$, $q\coloneqq u^4-v^5$, $R\coloneqq|x-y|$, and $S_m(a,b)\coloneqq\sum_{j=0}^{m-1}a^{m-1-j}b^j$. The pulled-back ball support function is represented by
\[
P\coloneqq\bigl(\bar qS_4(u,u'),-\bar qS_5(v,v')+\bar v^8S_8(v,v'),\bar w\bigr),
\quad \Phi\coloneqq P\cdot(x-y),\quad \Psi\coloneqq\Phi-\rho(x).
\]
These definitions give
\[
\Phi(x,y)=\sum_{j=1}^3\overline{F_j(x)}(F_j(x)-F_j(y)),\qquad
\Psi(x,y)=1-\sum_{j=1}^3\overline{F_j(x)}F_j(y).
\]
For $n=3$ and $h=1$, formula \cite[p.~125, (3.3)]{BP} expresses the solution as a sum of volume integrals with integrands of the following types:
\begin{equation}\label{eq:types}
\frac{\rho f\wedge\psi}{\Psi R^5},\qquad
\frac{f\wedge\db\rho\wedge\psi'}{\Psi^2R^3},\qquad
\frac{\rho f\wedge\psi''}{\Psi^3R^2}.
\end{equation}
Here their $r$ and $\widetilde\Phi$ are written as $\rho$ and $\Psi$. The coefficients of $\psi,\psi',\psi''$ are asserted to be smooth and bounded on $\Omega_3\times\Omega_3$.

We expand the volume integrand in the weighted representation \cite[Proposition~3.1]{BP}, omitting the common holomorphic volume form $du\wedge dv\wedge dw$ and normalization constant. It consists of a weighted Bochner--Martinelli term and a transition term. With $B\coloneqq(\bar x-\bar y)/R^2$ and $\db=\db_x$, integration of the homotopy between $B$ and $P/\Phi$ gives the transition form
\begin{equation}\label{eq:T}
T\coloneqq\frac{\det(B,P,\db P)}{2\Phi^2}
 +\frac{\det(B,P,\db B)}{2\Phi}.
\end{equation}
The term involving $\db(1/\Phi)$ drops out because it produces a column proportional to $P$. With $n=3$ and $h=1$, the weight in \cite[Proposition~3.1]{BP} is $a=1-(\Phi/\Psi)^6$, so
\[
\db a=6\Phi^5\Psi^{-7}(\rho\db\Phi-\Phi\db\rho).
\]
In this representation, the transition form contributes $f\wedge\db a\wedge T$ to the volume integrand.

For small $t>0$, set $W_t\coloneqq\sqrt{1-2t^{64}}$, $x_t\coloneqq(t^5,t^4,W_t)$, $y_t\coloneqq(0,0,W_t)$, and $f\coloneqq d\bar v$. Both points belong to $\Omega_3$, and
\[
\rho(x_t)=-t^{64},\quad \rho(y_t)=-2t^{64},\quad
\Phi=t^{64},\quad\Psi=2t^{64},\quad R^2=t^8(1+t^2).
\]
Along this path, $B=(t^5,t^4,0)/R^2$, $P=(0,t^{60},W_t)$, and differentiation before restriction gives $\partial_{\bar u}P=(4t^{30},-4t^{31},0)$. Consequently,
\[
[d\bar u]\det(B,P,\db P)=4W_tt^{26},\qquad
 f\wedge\db a=-\frac3{64}W_tt^{-64}d\bar v\wedge d\bar w.
\]
Here $[d\bar u]\omega$ denotes the $d\bar u$-coefficient of $\omega$. The second determinant in \eqref{eq:T} has $d\bar u$ coefficient $W_tt^{-12}/(1+t^2)^2$ and contributes $O(t^{-140})$ to the volume integrand, whereas the Bochner--Martinelli part is $O(R^{-5})=O(t^{-20})$, with $a=63/64$ along this path. Adding the transition and Bochner--Martinelli contributions, the coefficient of $d\bar u\wedge d\bar v\wedge d\bar w$ in the volume integrand is
\begin{equation}\label{eq:total}
-\frac3{32}W_t^2t^{-166}+O(t^{-140}).
\end{equation}
This nonzero leading term comes from $\det(B,P,\db P)/(2\Phi^2)$ in \eqref{eq:T}. With bounded coefficients, the three types in \eqref{eq:types} are respectively $O(t^{-20})$, $O(t^{-140})$, and $O(t^{-136})$ along the same path. Their sum is therefore $O(t^{-140})$ and cannot equal \eqref{eq:total}. Thus the volume integrand in \cite[Proposition~3.1]{BP} does not admit the direct bounded-coefficient decomposition claimed in the derivation of \cite[p.~125, (3.3)]{BP}.

\subsection{Incompatibility of the pointwise division estimates}\label{app:division}
The domain $\Omega_3$ and the smooth forms in Theorem~\ref{thm:main}(ii) contradict the sup-norm estimate asserted by Grundmeier--Simon--Stens\o nes~\cite[Main Theorem~1]{GSS}. An erroneous vanishing factor appears in their correction estimates. Using the notation of~\cite{GSS}, with $x=(\xi,z,w)$ and $2L$ the D'Angelo type, the weight in \cite[p.~562]{GSS} gives
\[
e^{\psi_1/2}
=e^{\widetilde\psi_0/2}\sqrt{|w|^2+\|x\|^{4L}}\,\|x\|,
\]
whereas the estimate for $v_1$ in \cite[p.~564]{GSS} uses $|w|\|x\|$ in its numerator. At $w=0$ and $x\ne0$, the factor $\sqrt{|w|^2+\|x\|^{4L}}\,\|x\|$ equals $\|x\|^{2L+1}>0$ and cannot supply the asserted vanishing in $w$.

That vanishing is also incompatible with the division construction itself. In \cite[pp.~539, 542, 555]{GSS}, the function
\[
h_1=\frac1\Phi-zv_1-wv_2
\]
is required to be holomorphic on an enlarged domain containing a punctured neighborhood of the base point in the original boundary. Test this on the ball $\operatorname{Re}\xi+|\xi|^2+|z|^2+|w|^2<0$, with base point $0$. The initial support function of \cite[p.~543]{GSS} is then $\Phi=\xi-A(|z|^2+|w|^2)$, $A>0$. Near a boundary point with $w=0$ and $z\ne0$, all denominator factors in \cite[Lemma~9.1]{GSS} are bounded below: here $k=L=1$, $M_1=|z|+\|x\|^2$, $M_2=|w|+\|x\|^2$, and $\operatorname{Re}\Phi=-|\xi|^2-(1+A)(|z|^2+|w|^2)<0$. Their estimates imply $v_1=O(|w|)$ and $v_2=O(1)$, forcing, on a sufficiently small relatively open neighborhood $S$ of such a point in $\{w=0,\ \operatorname{Re}\xi+|\xi|^2+|z|^2=0\}$ with $z\ne0$ and $1/2+\xi\ne0$,
\[
h_1(\xi,z,0)=\frac1{\xi-A|z|^2}
\quad\text{on }S.
\]
But the tangential CR field $\overline L\coloneqq\partial_{\bar z}-z(1/2+\xi)^{-1}\partial_{\bar\xi}$ annihilates the restriction of a holomorphic function to $S$, whereas
\[
\overline L\left(\frac1{\xi-A|z|^2}\right)
=\frac{Az}{(\xi-A|z|^2)^2}\ne0.
\]
Thus the claimed pointwise estimates cannot hold together with the required holomorphic division identity.

\section*{Funding}
The author is supported by the CAS Project for Young Scientists in Basic Research, Grant No. YSBR-033.

\section*{Competing interests}
The author declares no competing interests.

\section*{AI Use Disclosure}
GPT-6-Astra was used to find the counterexample, draft and revise the manuscript. The author has verified all mathematical statements and takes full responsibility for the entire content of the paper.


\begin{thebibliography}{29}
\bibitem{BP} Elisabetta Barletta and Carla Parrini, \emph{Bounded solutions for $\bar\partial$-problem in pseudo-Siegel domains}, Ann. Mat. Pura Appl. (4) \textbf{168} (1995), 119--132. \url{https://doi.org/10.1007/BF01759256}.
\bibitem{Berndtsson} Bo Berndtsson, \emph{A smooth pseudoconvex domain in $\C^2$ for which $L^\infty$-estimates for $\bar\partial$ do not hold}, Ark. Mat. \textbf{31} (1993), 209--218.
\bibitem{Catlin} David Catlin, \emph{Subelliptic estimates for the $\bar\partial$-Neumann problem on pseudoconvex domains}, Ann. of Math. (2) \textbf{126} (1987), no.~1, 131--191.
\bibitem{Cumenge} Anne Cumenge, \emph{Sharp estimates for $\bar\partial$ on convex domains of finite type}, Ark. Mat. \textbf{39} (2001), no.~1, 1--25. \url{https://doi.org/10.1007/BF02388789}.
\bibitem{DAngelo} John P. D'Angelo, \emph{Real hypersurfaces, orders of contact, and applications}, Ann. of Math. (2) \textbf{115} (1982), no.~3, 615--637.
\bibitem{DFischer} Klas Diederich and Bert Fischer, \emph{H\"older estimates on lineally convex domains of finite type}, Michigan Math. J. \textbf{54} (2006), no.~2, 341--352. \url{https://doi.org/10.1307/mmj/1156345598}.
\bibitem{DFF} Klas Diederich, Bert Fischer, and John Erik Forn\ae ss, \emph{H\"older estimates on convex domains of finite type}, Math. Z. \textbf{232} (1999), no.~1, 43--61.
\bibitem{DFRealAnalytic} Klas Diederich and John E. Forn\ae ss, \emph{Pseudoconvex domains with real-analytic boundary}, Ann. of Math. (2) \textbf{107} (1978), no.~3, 371--384. \url{https://doi.org/10.2307/1971120}.
\bibitem{DF} Klas Diederich and John Erik Forn\ae ss, \emph{Stein neighborhoods for finite preimages of regular domain}, Manuscripta Math. \textbf{50} (1985), 11--27. \url{https://doi.org/10.1007/BF01168825}.
\bibitem{DFSupport} Klas Diederich and John Erik Forn\ae ss, \emph{Support functions for convex domains of finite type}, Math. Z. \textbf{230} (1999), no.~1, 145--164. \url{https://doi.org/10.1007/PL00004683}.
\bibitem{DFW} Klas Diederich, John Erik Forn\ae ss, and Jan Wiegerinck, \emph{Sharp H\"older estimates for $\bar\partial$ on ellipsoids}, Manuscripta Math. \textbf{56} (1986), no.~4, 399--417.
\bibitem{FeffermanKohn} Charles L. Fefferman and Joseph J. Kohn, \emph{H\"older estimates on domains of complex dimension two and on three-dimensional CR manifolds}, Adv. in Math. \textbf{69} (1988), no.~2, 223--303.
\bibitem{FKM} Charles L. Fefferman, Joseph J. Kohn, and Matei Machedon, \emph{H\"older estimates on CR manifolds with a diagonalizable Levi form}, Adv. Math. \textbf{84} (1990), no.~1, 1--90. \url{https://doi.org/10.1016/0001-8708(90)90036-M}.
\bibitem{Fornaess} John Erik Forn\ae ss, \emph{Sup-norm estimates for $\bar\partial$ in $\C^2$}, Ann. of Math. (2) \textbf{123} (1986), no.~2, 335--345.
\bibitem{FornaessKim} John-Erik Fornaess and Kang-Tae Kim, \emph{Some problems}, Complex Analysis and Geometry: KSCV10, Gyeongju, Korea, August 2014 (Filippo Bracci, Jisoo Byun, Herv\'e Gaussier, Kengo Hirachi, Kang-Tae Kim, and Nikolay Shcherbina, eds.), Springer Proceedings in Mathematics \& Statistics, vol.~144, Springer, Tokyo, 2015, pp.~369--377. \url{https://doi.org/10.1007/978-4-431-55744-9_28}.
\bibitem{FS} John Erik Forn\ae ss and Nessim Sibony, \emph{Smooth pseudoconvex domains in $\C^2$ for which the corona theorem and $L^p$ estimates for $\bar\partial$ fail}, Complex Analysis and Geometry (Vincenzo Ancona and Alessandro Silva, eds.), The University Series in Mathematics, Plenum Press, New York, 1993, pp.~209--222.
\bibitem{Gavosto} Estela A. Gavosto, \emph{H\"older estimates for the $\bar\partial$-equation in some domains of finite type}, J. Geom. Anal. \textbf{7} (1997), no.~4, 593--609. \url{https://doi.org/10.1007/BF02921635}.
\bibitem{GrauertLieb} Hans Grauert and Ingo Lieb, \emph{Das Ramirezsche Integral und die L\"osung der Gleichung $\bar\partial f=\alpha$ im Bereich der beschr\"ankten Formen}, Rice Univ. Stud. \textbf{56} (1970), no.~2, 29--50.
\bibitem{GSS} Dusty Grundmeier, Lars Simon, and Berit Stens\o nes, \emph{Sup-norm estimates for $\bar\partial$}, Pure Appl. Math. Q. \textbf{18} (2022), no.~2, 531--571. \url{https://doi.org/10.4310/pamq.2022.v18.n2.a8}.
\bibitem{Henkin} G. M. Henkin, \emph{Integral representation of functions in strictly pseudoconvex domains and applications to the $\bar\partial$-problem}, Math. USSR-Sb. \textbf{11} (1970), no.~2, 273--281, English translation of Mat. Sb. (N.S.) \textbf{82(124)} (1970), 300--308.
\bibitem{Kerzman} Norberto Kerzman, \emph{H\"older and $L^p$ estimates for solutions of $\bar\partial u=f$ in strongly pseudoconvex domains}, Comm. Pure Appl. Math. \textbf{24} (1971), no.~3, 301--379.
\bibitem{Kohn} Joseph J. Kohn, \emph{Boundary behavior of $\bar\partial$ on weakly pseudo-convex manifolds of dimension two}, J. Differential Geom. \textbf{6} (1972), no.~4, 523--542.
\bibitem{McNeal} Jeffery D. McNeal, \emph{Convex domains of finite type}, J. Funct. Anal. \textbf{108} (1992), no.~2, 361--373.
\bibitem{PSP} Carla Parrini and Angela Selvaggi Primicerio, \emph{Uniform estimate for the $\bar\partial$-problem in a class of branched covering of strictly pseudoconvex domains}, Math. Z. \textbf{189} (1985), 465--474. \url{https://doi.org/10.1007/BF01168153}.
\bibitem{Ramirez} Enrique Ram\'irez de Arellano, \emph{Ein Divisionsproblem und Randintegraldarstellungen in der komplexen Analysis}, Math. Ann. \textbf{184} (1970), 172--187.
\bibitem{Range} R. Michael Range, \emph{Integral kernels and H\"older estimates for $\bar\partial$ on pseudoconvex domains of finite type in $\C^2$}, Math. Ann. \textbf{288} (1990), no.~1, 63--74.
\bibitem{Sibony} Nessim Sibony, \emph{Un exemple de domaine pseudoconvexe r\'egulier o\`u l'\'equation $\bar\partial u=f$ n'admet pas de solution born\'ee pour $f$ born\'ee}, Invent. Math. \textbf{62} (1980), no.~2, 235--242.
\bibitem{Yao} Liding Yao, \emph{Sobolev and H\"older estimates for homotopy operators of the $\bar\partial$-equation on convex domains of finite multitype}, J. Math. Anal. Appl. \textbf{538} (2024), no.~2, 128238. \url{https://doi.org/10.1016/j.jmaa.2024.128238}.
\end{thebibliography}
\end{document}